\documentclass[12pt]{article}
\usepackage{amssymb}
\def\hang{\hangindent\parindent}
\def\tex#1{\indent\llap{[#1]\enspace}\ignorespaces}
\def\re{\par\hang\tex}

\def\d{\delta}
\def\D{\Delta}
\def\UU{{\cal U}}
\def\WW{{\cal W}}

\def\Der{{\rm Der}}

\def\p{\psi}
\def\Ker{{\rm Ker}}
\def\es{\varepsilon}
\def\Im{{\rm Im}}
\def\span{{\rm span}}

\def\v{\varphi}

\def\ssc{\scriptscriptstyle}

\def\cl{\centerline}

\def\ul{\underline}

\def\rar{\rightarrow}

\def\bs{\backslash}

\def\vs{\vspace*}

\def\ra{\rangle}

\def\la{\langle}
\def\ni{\noindent}
\def\ptl{\partial}

\def\Z{\mathbb{Z}{\ssc\,}}

\def\F{\mathbb{F}{\ssc\,}}
\begin{document}
%
%\setcounter{page}{0}
%$[x,a_i]\otimes b_i$ need be changed to $\v(x,a_i)\otimes b_i$
%after Definition 2.3.
%
%Prop 3.5 need to be changed: $r$ should in $\WW\otimes\WW\otimes\WW$.
%\pagebreak
%
%
\cl{{\Large \bf
 Lie bialgebras of generalized Witt type}\footnote{Supported by NSF grant no.~10471091 of China, two grants ``Excellent Young Teacher Program'' and ``Trans-Century Training Programme Foundation for the Talents'' from Ministry of Education of China}}
\vs{6pt}

\cl{ Guang'ai Song and Yucai Su}
\cl{\small Department of Mathematics, Shanghai Jiao Tong University,
 Shanghai, 200240, \vs{-4pt}China}
\cl{\small E-mail: gasong@sjtu.edu.cn, ycsu@sjtu.edu.cn}
\vs{6pt}

{\small
\parskip .005 truein
\baselineskip 3pt \lineskip 3pt

\noindent{{\bf Abstract.} In a paper by Michaelis a class
of infinite-dimensional Lie bialgebras containing the Virasoro
algebra was presented.
This type of Lie bialgebras was classified by Ng and Taft.
In this paper, all Lie bialgebra structures on
the Lie algebras of generalized Witt type are classified.
It is proved that,
for any Lie algebra $\WW$
of generalized Witt type,
all Lie bialgebras on $\WW$
are coboundary triangular Lie bialgebras. As a by-product, it is
also proved that
the first cohomology group
$H^{1}(\WW,\WW \otimes\WW)$ is trivial.

\vs{5pt}

}

}

\vs{6pt}

\cl{\bf\S1. \
Introduction}\setcounter{section}{1}\setcounter{equation}{0}
Recently there appeared a number of papers on the structure theory
of Lie bialgebras. In general, a Lie bialgebra is a vector space
%which is
endowed simultaneously with the structure of a Lie algebra
and the structure of a Lie coalgebra, and some compatibility condition, which was suggested by a study of Hamiltonian mechanics and Poisson
Lie group, holds.
Michaelis [M] %[M1] 
presented a class of Witt type Lie
bialgebras. Michaelis also gave a method on how
to obtain the structure of a triangular, coboundary Lie bialgebra
on a Lie algebra which contains two linear independent elements $a$ and
$b$ satisfying condition $[a,b]=kb$ for some non-zero
scalar $k$.
This kind of Lie bialgebras was classified by Ng and Taft [NT]
(cf.~[N], [T]).
In this paper, we classified all
Lie bialgebra structures on the generalized Witt Lie algebras which
are defined in [DZ] (see also [P], [SXZ], [SZ], [X]).
 \vskip10pt
\cl{\bf\S2. \ Preliminaries}\setcounter{section}{2} We start with
recalling the some concepts of Lie bialgebras. Throughout the
paper, $\F$ denotes a field of characteristic zero. Let $L$ be a
vector space over $\F$. Denote by $\tau$ the twist map of $L
\otimes L$, namely\vs{-4pt},
$$
\tau(x\otimes y)= y \otimes x\mbox{ \ \  for \ }x,y\in L\vs{-4pt}.
$$
Denote by $\xi$ the map which cyclically permutes the coordinates
of $L \otimes L \otimes L$, i.e.\vs{-4pt},
$$
\xi (x_{1} \otimes x_{2}
\otimes x_{3}) =x_{2} \otimes x_{3} \otimes x_{1}
\mbox{ \ \ for \ }x_1,x_2,x_3\in L.
$$
\vskip0pt

\ni{\bf Definition 2.1.} Let $L$ be a vector space over $\F$ and
let $\v :L \otimes L \rar L$ be a bilinear map. The pair $(L, \v)$
is called a {\it Lie algebra} if the following conditions are
satisfied: \vs{-3pt}\begin{itemize}\parskip-3pt
\item[(i)]
$\Ker(1-\tau) \subset \Ker\,\v$ (where $1$ is the identity map
of $L \otimes L$);
\item[(ii)]
$\v \cdot (1 \otimes \v ) \cdot (1 + \xi +\xi^{2}) =0 : L
\otimes L \otimes L \rar L.$
\end{itemize}\vs{-3pt}

Note that $\Ker(1- \tau) = \span\{x \otimes x\, |\,x \in L\}$. Also
$\Im(1 + \tau) \subset \Ker(1 - \tau)$ and the equality holds if
the characteristic of $\F$ is not $2$. In this case, we can replace condition (i) of Definition 2.1 by condition

{(i$^{\prime} ): \ \v = -\v \cdot\tau.$} \vskip5pt

\ni {\bf Definition 2.2.} Let $M$ be a vector space over $\F$ and
let $\D: M \rar M \otimes M$ be a linear map. The pair $(M, \D)$
is called a {\it Lie coalgebra} if the following conditions are
satiafied: \vs{-3pt}\begin{itemize}\parskip-3pt
\item[(i)] $\Im\,\D \subset \Im(1- \tau)$,
\item[(ii)]
$(1 + \xi +\xi^{2}) \cdot (1 \otimes \D) \cdot \D =0: M \rar M
\otimes M \otimes M.$
\end{itemize}\vs{-3pt}

The map $\D$ is called the {\it comultiplication} or
{\it cobracket} or {\it diagonal} of $M$.
Condition (i)  of Definition 2.2 is called the {\it strong anti-commutativity},
and condition (ii) is called the {\it Jacobi-identity}.
Similar to Definition 2.1,
since $\Im(1-\tau) \subset \Ker(1 + \tau)$,
 and the equality holds if
the characteristic of $\F$ is not $2$, we can replace condition (i)
 of Definition 2.2 by condition

{(i$^{\prime} ): \ \D =-\tau \D.$} \vskip5pt

\ni {\bf Definition 2.3.} A {\it Lie bialgebra} is a triple $(L,
\v, \D )$ which satisfies the following conditions
\vs{-3pt}\begin{itemize}\parskip-3pt
\item[(i)] $(L, \v)$ is a Lie algebra,
\item[(ii)] $(L, \D)$ is a Lie coalgebra,
\item[(iii)]
$\D  \v (x, y) = x \cdot \D y - y \cdot \D x$ for $x, y \in
L$,\vs{-3pt}
\end{itemize}
\noindent where the symbol ``$\cdot$'' means the action
$x\cdot (\sum_{i} {a_{i} \otimes b_{i}}) = \sum_{i} ( {[x,
a_{i}] \otimes b_{i} + a_{i} \otimes [x, b_{i}]})$
for $x, a_{i}, b_{i} \in L$, and
in general $[x,y]=\v(x,y)$ for $x,y
\in L$.
\vskip7pt

\ni
{\bf Remark 2.4.} The difference between definition of a
Lie bialgebra and definition of a bialgebra lies in
the compatibility condition (iii): The compatibility condition
of a bialgebra requires that
 $\D$ is an algebra
morphism, i.e. $\D \cdot \v = (\v \otimes \v) \cdot (1 \otimes
\tau \otimes 1) \cdot \D \otimes \D$, however
the compatibility condition of a Lie bialgebra requires that
$\D$ is a derivation of $L \rar L
\otimes L.$ Thus the properties of Lie bialgebras are not similar
to those of bialgebras.
\vskip7pt

\ni
{\bf Definition 2.5.} A {\it coboundary Lie bialgebra}
is a $(L, \v, \D,r),$
where $(L, \v, \D)$ is a Lie bialgebra and $r \in  \Im(1 -
\tau) \subset L \otimes L$ such that $\D$ is a {\it coboundary of $r$},
i.e. for arbitrary $x \in L$,
$$\D (x) = x \cdot r.$$
%\vskip3pt

\ni
{\bf Definition 2.6.} A coboundary Lie bialgebra $(L, \v,\D, r)$
is called {\it triangular} if it satisfies the
following {\it classical Yang-Baxter
Equation} (CYBE): $$c(r)=0,$$
 where $c(r)$ is defined by
\begin{equation}
\label{add1-}
\mbox{$c(r): = [r^{12} , r^{13}] +[r^{12} , r^{23}] +[r^{13} ,
r^{23}],$}
\end{equation}
and $r^{ij}$ are defined as follows:
Denote by $\UU(L)$ the {\it universal enveloping algebra} of
$L$. If $r =\sum_{i} {a_{i} \otimes b_{i}} \in
L \otimes L$, then
$$\begin{array}{c}
 r^{12} = \sum \limits_{i}{a_{i} \otimes b_{i}
\otimes 1} \in {\UU (L) \otimes \UU (L) \otimes \UU (L)}, \\
r^{13}= \sum \limits_{i} {a_{i} \otimes 1 \otimes b_{i}} \in{ \UU
(L) \otimes \UU (L) \otimes \UU (L)}, \\
r^{23} = \sum \limits_{i}{1 \otimes a_{i} \otimes b_{i}} \in \UU
(L) \otimes \UU (L) \otimes \UU(L).
\end{array}$$
Here $1$ is the identity element of $\UU (L)$. Obviously,
$$
\begin{array}{l}
[r^{12} ,r^{13}] = \sum\limits_{i, j} [a_{i}, a_{j}] \otimes
b_{i} \otimes b_{j}\in L \otimes L \otimes L,\\[4pt]
[r^{12}
,r^{23}] = \sum\limits_{i,k} a_{i} \otimes[ b_{i},a_{k}] \otimes
b_{k}\in L \otimes L \otimes L,\\[4pt]
[r^{13} ,r^{23}] = \sum\limits_{j,k} a_{j} \otimes a_{k}
\otimes[ b_{j}, b_{k}] \in L \otimes L \otimes L.
\end{array}$$
\vskip3pt

Let $L$ be a Lie algebra, then $L \otimes L$ is an $L$-module under the adjoint
diagonal action of $L$.
If $r \in \Im(1-\tau) \subset L \otimes
L$, we define a linear map
$$\D = \D_{r} : L \rar L \otimes L$$
via
\begin{equation}
\label{add2.0}
\D_r (x) = x \cdot r.
\end{equation}
\ni
Then it is obvious that $\Im\D \subset\Im(1 - \tau).$ Furthermore,
$$\D [x, y] =x \cdot \D y - y \cdot \D x.$$
Thus, $(L, [\cdot ,\cdot ], \D)$ is a
coboundary Lie bialgebra if $\D$ satisfies Jacobi identity.

The following result is due to Drinfeld [D].
\vskip7pt

\ni{\bf Theorem 2.7.}
{\it Let $L$ be a Lie algebra then $\D = \D_{r}$ $($for some $r
\in\Im(1-\tau))$ endows $(L,[\cdot,\cdot], \D )$ with a Lie bialgebra
structure if and only if $r$ satisfies the following modern Yang-Baxter
Equation {\rm(MYBE):}
$$x \cdot c(r) =0\mbox{ \ \ \ for all \ \ }x\in L.
\eqno(2.4)$$
}
%\vskip-5pt

The following main result of [M] %[M1] 
is a slight generalization of Theorem 2.7
under a rather strong condition.
\vskip7pt

\ni
{\bf Proposition 2.8.}  {\it Let $L$ be a Lie algebra containing two
linear
independent elements $a$ and $b$ satisfying $[a, b] = kb$ with $0 \neq
k \in \F$.
Set
$$r = a \otimes b -b \otimes a$$
and define a linear map
$$\D_{r} (x) = x \cdot r = [x, a] \otimes b - b
\otimes [x, a] + a \otimes [x, b] - [x, b] \otimes a
\mbox{ \ \ for \ }x \in L.$$
Then $\D_{r}$ equips $L$ with the
structure of a triangular coboundary Lie bialgebra.
}\vskip7pt

The following result was obtained in [NT].
\vskip7pt

\setcounter{equation}{4} \ni {\bf Proposition 2.9.} {\it Let $L$
be a Lie algebra. Set $r =\sum_{i=1} ^{n}( a_{i} \otimes b_{i} -
b_{i} \otimes a_{i}) \in L \otimes L$ for any $a_i,b_i\in L$, and
set $\D = \D_{r}$. Then for any $x \in L$,
\begin{equation}
\label{add-c}
(1 + \xi + \xi^{2}) \cdot (1 \otimes \D) \cdot \D (x) = x
\cdot c (r),
\end{equation}
where $c(r)$ is defined by $(2.1)$. 
In particular, 
%Consequently 
for any $x \in L$,
$$(1 + \xi + \xi^{2}) \cdot (1 \otimes \D) \cdot \D (x) = 0
\mbox{ \ \ if
and only if  \ }x \cdot c (r) =0.$$
%In particular if $c (r) =0 $ then $\D$ satisfies the Jacobi-identity.
}
%\vskip15pt

\cl{\bf\S3. \  The structures of Lie bialgebras of generalized Witt
type}\setcounter{section}{3}\setcounter{equation}{0}
Let $\F$ be a field of characteristic zero, and let $n>0$.
Let $A$ be a {\it nondegenerate} additive
subgroup of $\F^n$, i.e., $A$ contains an $\F$-basis
%$\{\e_{1}, \e_{2}, \cdots, \e_{n}\}$
of the vector space $\F^n$.
%Fix $\g = a_{1} \e_{1}+a_2\e_2+\cdots+a_{n} \e_{n} \in A$
% with $a_{i} \neq 0$ for all $i = 1,2,\cdots ,n.$

Let $\F[A] =\span\{ t^{x}\, |\, x \in A \}$ be the group algebra of $A$ over $\F$ such that
$t^x\cdot t^y=t^{x+y}$ for $x,y\in A$.
Let $T=\span\{\ptl_i\,|\,i=1,2,...,n\}$ be an $n$-dimensional vector space over $\F$.
Regarding elements of $T$ as derivations of $\F[A]$ by setting
$$\ptl_{i} t^{x} = x_{i} t^{x}
\mbox{ \ for \ }x=(x_1,x_2,...,x_n)\in A\subset\F^n.$$
 If $A = \Z^{n}$, then
$\F[A]$ agrees with the Laurent polynomial algebra $\F [t_{1} ^{\pm 1},
\cdots, t_{n} ^{\pm 1} ]$, and $\ptl_{i}$ coincides with the {\it degree operator}
$t_{i}\frac {\ptl}{\ptl_{t_{i}}}$.

Denote $\WW = \F[A] \otimes T =\span\{t^{x} \ptl \,|\,
x\in A, \ptl \in T \}$, where we have simplified the notation by setting
$t^{x} \ptl=t^{x}\otimes \ptl$. Then $\WW$ is a simple
{\it Lie algebra of
generalized Witt type} [DZ] (see also [SXZ], [SZ], [X]) under the following bracket
$$
\mbox{$[t^{x} \ptl , t^{y} \ptl^{\prime}] = t^{x + y} (\ptl (y)
\ptl^{\prime} - \ptl^{\prime} (x) \ptl)$ \ \ for \ $x,y\in A,\ \ptl,\ptl'\in T,$}
\eqno(3.1)$$
where in general,
$$
\ptl(y)=\la\ptl,y\ra=\la y,\ptl\ra=\mbox{$\sum\limits_{i=1}^n$}a_iy_i
\mbox{ \ \ for \ $\ptl=\sum\limits_{i=1}^n a_i\ptl_i\in T,\,y=(y_1,y_2,...,y_n)\in A$}.
\eqno(3.2)$$
There is a natural
$A$-gradation on $\WW=\oplus_{x\in A}\WW_x$
by setting $\WW_{x} = \{t^{x} \ptl\, |\, \ptl \in
T \}$ for $x \in A$. This gradation is also compatible with
the Lie algebra structure, i.e., $[\WW_{x} , \WW_{y}] \subset
\WW_{x + y}$, and $\WW_{0} = T$ is a torus.
The decomposition $\WW=\oplus_{x\in A}\WW_x$
is also the root space decomposition with respect to the torus $\WW_{0}$.
Obviously, $\WW$ is not necessarily finitely-generated.

\setcounter{equation}{2}
The bilinear map $\la\cdot,\cdot\ra:T \times A \rar F$ in (3.2)
is {\it nondegenerate} in the sense
\begin{eqnarray}
\label{3.1}
\!\!\!\!\!\!\!\!\!\!\!\!&\!\!\!\!\!\!\!\!\!&
\ptl(x) = \la\ptl, x\ra =0 \mbox{ \ for any \ } \ptl \in T \ \ \Longrightarrow\ \
x=0,\mbox{ \ and}
\\
\label{3.2}
\!\!\!\!\!\!\!\!\!\!\!\!&\!\!\!\!\!\!\!\!\!&
\ptl(y) = \la\ptl, y\ra = 0  \mbox{ \ for any \ } y \in A\ \ \Longrightarrow\ \
\ptl=0.
\end{eqnarray}

The tensor product $V=\WW \otimes \WW$ is an $A$-graded $\WW$-module
under the adjoint diagonal action of $\WW$. The gradation is given by
\begin{equation}
\label{3.3}
V= \bigoplus_{x \in A} V_{x},
\mbox{ \ where \
$V_{x} = \sum \limits_{x_{1} +x_{2} =x} \WW_{x_{1}} \otimes\WW_{x_{2}}$.}
\end{equation}

We shall discuss the derivation algebra $\Der(\WW, V).$ First let us
recall some basic definitions.
\vskip7pt

\ni{\bf Definition 3.1.}
Let $L=\oplus_{x\in A}L_x$ be an $A$-graded Lie algebra for some abelian group
$A$, and let $V=\oplus_{x\in A}V_x$ be an $A$-graded $L$-module.
A linear map $D: L \rar V$ is called a {\it derivation} if it satisfies

\cl{$D([u, w]) = u \cdot D (w) - w \cdot D(u)$ \ for $u, w \in L.$}

\noindent
A derivation of the form

\cl{$ D (u) = u \cdot v$ \ for  $u \in L$ (with a fixed $v \in V$),}

\noindent  is called an {\it inner
derivation}. A derivation $D$ is  {\it homogeneous of degree
$x$} if $D (V_{y}) \subset V_{x +y}$ for all $y \in A$.
Denote by $\Der(L, V)$ and by ${\rm Inn}(L, V)$ the
space of derivations and the space of inner derivations respectively.
Denote by $\Der(L, V)_{x} = \{D \in \Der(L, V) \,|\,{\rm deg\,} D =x\}$.
\vskip5pt

It is well known that the first cohomology group of $L$ with coefficients in the module
$V$ is isomorphic to
$$
H^{1} (L, V) \cong \Der(L, V)/{\rm Inn}(L, V).
$$

The following two propositions
are slight generalizations of some results in [F].
First, by dropping the condition of $L$ being finitely-generated,
we can generalize
[F, Proposition 1.1] to obtain
\vskip7pt

\ni{\bf Proposition 3.2.}
{\it Let $L$ be a Lie algebra, and let
$V$ be an $A$-graded $L$-module. Then
$$
\mbox{$d = \sum\limits_{x \in A} d_x , \mbox{  \ where \ }d_x \in \Der(L,
V)_x ,$}
$$
which holds in the sense that for every $u \in L$, only finitely
many $d_x (u)\neq 0,$ and $d(u) = \sum_{x \in A} d_x(u). $
}\vskip7pt

\ni{\bf Proposition 3.3.}
{\it Let $L$ be a Lie algebra, and let $V$ be an
$A$-graded $L$-module such that
\begin{itemize}\parskip-3pt
\item[{\rm(a)}] $H^{1} (L_{0} , V_{x}) = 0$ for $x \in A \setminus \{ 0\},$
\item[{\rm(b)}] ${\rm Hom}_{L_{0}} (L_{x} , V_{y}) = 0$ for $x \neq y$,
\item[{\rm(c)}] ${\rm dim\,}L_{0} = n < \infty.$
\end{itemize}
\noindent Then

\cl{$\Der(L, V)=\Der(L, V)_{0} +{\rm Inn} (L, V).$}
}\vskip5pt

\ni{\it Proof.}
Let $d \in \Der(L, V).$ According to Proposition 3.2, we
have
$$
\mbox{$d = \sum\limits_{x \in A} d_x.$}
$$
Suppose $x \neq 0.$ Then $d_x |_{L_{0}}$ is a derivation
from $L_{0}$ to the $L_{0}$-module $V_x.$ By virtue of (a),
$d_x |_{L_{0}}$ is inner, i.e., there exists $v_x \in V_x$ such
that $d_x (u) = u \cdot v_x$ for $u \in L_0.$ Let $u_{1} , u_2,\cdots,
u_{n}$ be a basis of $L_{0}$. Let $S \subset A$ be the set of all elements
$x \neq 0$ such that $v_x \neq 0,$ and for $i=1,2,...,n$, let $S_{i}$
be the set of all elements $x \in S$ such that $u_{i} \cdot v_x \neq
0.$ For $x \in S, $ we have

 \cl{$d_x (u_{i}) = u_{i} \cdot v_x .$}

\noindent By Proposition 3.2,
 $d_x (u_{i}) \neq 0$ for only
finite many $x \in A$, i.e., $S_{i}$ is a finite set. Thus
the union set
$S=\cup_{i=1}^n S_i$
is also finite.

If $x \neq 0$ with $d_x |_{L_{0}} = 0$, which means that $d_x \in
{\rm Hom}_{L_{0}} (L_{0} , V_x).$ From (b), we obtain $d_x = 0.$
This proves that
$$
x\ne0,\ d_x\ne0\ \ \Longrightarrow\ \ x\in S.
$$
Thus
$d=\sum_{x\in A} d_x=\sum_{x\in S\cup\{0\}}d_x$ is a finite sum.

It remains to prove that for $x \in S,$ $d_x$ is inner.
Consider $\p_x : L \rar V$ with $\p_x (u) = d_x (u) - u\cdot
v_x$ for $u \in L.$ Obviously, $\p_x$ is a homogeneous derivation of degree $x$
which vanishes on $L_{0},$ hence it is a homomorphism of
$L_{0}$-modules. From (b), we obtain that $\p_x = 0,$ namely, $d_x \in
{\rm Inn}(L, V)_x.$ \hfill $\Box$
\vskip7pt

\ni{\bf Proposition 3.4.}
{\it
Let $\WW$ be the Lie algebra of generalized Witt type defined in $(3.1)$. Let
$V = \WW \otimes \WW = \oplus_{x \in A} V_{x}$. Then
\begin{equation}
\label{3.4}
\Der(\WW, V) ={\rm Inn}(\WW, V) + \Der(\WW, V)_{0}
\end{equation}
}\vskip5pt\ni
{\it Proof.} From Proposition 3.3, we just need to verify the
conditions of Proposition 3.3. Condition (c) is obvious.

{\bf Step 1}: {\it First we verify condition {\rm(a)}.}

Suppose $D \in\Der (\WW_{0} , V_{x} )$ and $x \neq 0. $
Let
$\ptl^{\prime} , \ptl \in T =\WW_{0}.$ Since $[\ptl^{\prime} ,
\ptl] = 0,$ applying $D$ to it, we obtain
\begin{equation}
\label{3.5}
\ptl^{\prime}\cdot D (\ptl) - \ptl\cdot D (\ptl^{\prime}) = 0,
\end{equation}
where $D(\ptl), D(\ptl^{\prime}) \in V_{x}=\sum_{x_1+x_2=x}\WW_{x_1}\otimes\WW_{x_2}.$
Since $\ptl\cdot a=\ptl(x)a$ for all $a\in V_x$, from (\ref{3.5}) we obtain
\begin{equation}
\label{a1}
\ptl^{\prime} (x) D(\ptl) = \ptl(x) D(\ptl^{\prime})
\end{equation}
By (\ref{3.1}), we can choose $\ptl^{\prime} \in T$ such that
$\ptl^{\prime} (x) \neq 0.$ Set $w = (\ptl^{\prime} (x))^{-1}
D(\ptl^{\prime})\in V_x$. Then (\ref{a1}) shows
that $D(\ptl)=\ptl \cdot w$ for all $\ptl\in\WW_0$,
i.e., $D$ is inner. Thus

\cl{H$^{1} (\WW_{0} ,V_{x}) = 0. $}

{\bf Step 2}: {\it Next we verify condition {\rm(b)}.}

Let $x\ne y$, and let $f\in{\rm Hom}_{\WW_{0}} (\WW_{x} , V_{y})$. Then
\begin{equation}
\label{3.5+}
f (\ptl\cdot a) = \ptl \cdot f (a)\mbox{ \ \ for all \ }a\in\WW_x\mbox{ \ and \ }\ptl\in\WW_0.
\end{equation}
By (\ref{3.1}), there exists $\ptl \in \WW_{0}$ such that $\ptl
(x) \neq \ptl (y).$ The left-hand side of (\ref{3.5+}) is $\ptl
(x) f (a)$, but the right-hand side of (\ref{3.5+}) is $\ptl (y) f
(a)$ since $f(a)\in V_y$. This shows that $f = 0.$ \hfill $\Box$
\vskip7pt \ni 
{\bf Proposition 3.5.} \ 
{\it Suppose $c\in\WW\otimes\WW\otimes\WW$
satisfying $a\cdot c=0$ for all $a\in\WW$. Then $c=0$.}
\vskip5pt\ni 
{\it Proof.} 
Write $c=\sum_{x\in A}c_x$ as a finite
sum with $c_x\in (\WW\otimes\WW\otimes\WW)_x$. From $0=\ptl\cdot c=\sum_{x\in
A}\ptl(x)c_x$ for any $\ptl\in T$, we obtain $c=c_0\in V_0$. Now
write 
$$\mbox{
$c=\sum\limits_{x,y\in A}t^x\ptl_x\otimes t^y\ptl'_y\otimes t^{-x-y}\ptl''_{x,y}$ for some
$\ptl_x,\ptl'_y, \ptl''_{x,y}\in T$,
}$$
 where $\{(x,y)\in A\times A\,|\,
\ptl_x,\ptl'_y,\ptl''_{x,y}\ne0\}$
is a finite set. Choose a total order on $A$ compatible with its
group strcuture. 
Define the total order on $A\times A$ by: 
$$\mbox{
$(x_1,y_1)>(x_2,y_2)\ \Longleftrightarrow\
x_1>x_2$, or $x_1=x_2,\,y_1>y_2$.
}
$$
If $c\ne0$, let $(x_0,y_0)$ be the maximal element with
$\ptl_{x_0},\ptl'_{y_0},\ptl''_{x_0,y_0}\ne0$. In this case, the term $t^{x_0}\ptl_{x_0}
\otimes t^{y_0}\ptl'_{y_0}\otimes
t^{-x_0-y_0}\ptl'_{x_0,y_0}$ is called the {\it leading term} of $c$. Choose any
$z>0$ such that $\ptl_{x_0}(z-x_0)\ne0$. Then
$\ptl_{x_0}(z-x_0)t^{z+x_0}\ptl_{x_0}\otimes t^{y_0}\ptl'_{y_0}\otimes
t^{-x_0-y_0}\ptl'_{x_0,y_0}$ is the leading
term of $(t^z\ptl_{x_0})\cdot c$, a contradiction with the fact that
$(t^z\ptl_{x_0})\cdot c=0.$ \hfill$\Box$ \vskip7pt

\ni
{\bf Theorem 3.6.} \ {\it $\Der(\WW, V) ={\rm Inn}(\WW, V).$}
\vskip5pt\ni
{\it Proof.}
By Proposition 3.4, it suffices to prove that a derivation $D \in \Der (\WW, V)_{0}$ is inner.
We shall prove that after a number of steps in each
of which $d$ is replaced by $D-D'$ for some $D'\in {\rm Inn}(\WW,V)$ the 0 derivation is
obtained and thus proving that $D\in {\rm Inn}(\WW,V)$. This will be done by three claims.
\par
First using the fact that $A$ is a nondegenerate subgroup of $\F^n$, we can choose an
$\F$-basis $\{\es_1,\es_2,...,\es_n\}\subset A$ of the vector space $\F^n$.
We re-denote $\{\ptl_1,\ptl_2,...,\ptl_n\}\subset T$ to be the {\it dual basis} of $\{\es_1,\es_2,...,\es_n\}$ under the pairing (3.2), i.e.,
$$
\la\ptl_i,\es_j\ra=\d_{i,j}\mbox{ \ \ for \ \ }i,j=1,2,...,n.
\eqno(3.10)
$$
From now on, we shall denote an element $x$ of $A$
by
$$x=x_1\es_1+x_2\es_2+\cdots+x_n\es_n.$$
Recall that we have Lie bracket (3.1) for $\WW$.\vskip4pt

{\bf Claim 1}. $D(\ptl)=0$ for $\ptl\in T$.
\par
To prove this, applying $D$ to $[\ptl,t^x\ptl']=\ptl(x)t^x\ptl'$
for $x\in A,\,\ptl'\in T$, we obtain that \mbox{$(t^x\ptl')\cdot D(\ptl)=0$.}
By Proposition 3.5, $D(\ptl)=0$.
\vskip4pt

{\bf Claim 2}. We can suppose $D(t^{p{\ssc\,}\es_j}\ptl_j)=0$
for $j=1,2,...,n$ and $p=0,\pm1,-2$.
\par
We prove this by induction on $j$. Assume that $j$ is fixed and
suppose we have proved
$$D(t^{p{\ssc\,}\es_i}\ptl_i)=0 \mbox{ \ \ for \ \ }i<j,\ p=0,\pm1,
-2.
\eqno(3.11)$$
For any $u\in\WW_z$ with $z\in A$, we assume
$$
D(u)=\sum_{x\in A,\,
1\le k,\ell\le n}c^{(u)}_{x,k,\ell\,}t^{x+z}\ptl_k\otimes t^{-x}\ptl_{\ell\,}
\mbox{ \ \ for some \ }c^{(u)}_{x,k,\ell}\in\F,
\eqno(3.12)
$$
where $M_u=\{(x,k,\ell)\,|\,c^{(u)}_{x,k,\ell}\ne0\}$ is a finite set.
For $p\in\Z$,
we simply denote the coefficient
$c^{(u)}_{x,k,\ell}$ by $c^{(p)}_{x,k,\ell}$ when $u=t^{p{\ssc\,}\es_j}\ptl_j$.

Applying $D$ to $[t^{p{\ssc\,}\es_i}\ptl_i,t^{\es_j}\ptl_j]=0$ for $i<j$, by (3.11),
we obtain that if $c^{(1)}_{x,k,\ell}\ne0$ then
$$
x_i-p\d_{i,k}=x_i+p\d_{i,\ell}=0\mbox{ \ \ for \ \ }i<j,\ \,p=0,\pm1,-2.
\eqno(3.13)$$
For any $u\in V$, we denote by $u_d$
the inner derivation $w\mapsto w\cdot u,\,\forall\,w\in\WW$,
determined by $u$. Then for $x\in A$, we have
$$
(t^x\ptl_k\otimes t^{-x}\ptl_\ell)_d(t^{\es_j}\ptl_j)=
(x_j-\d_{j,k})t^{x+\es_j}\ptl_k\otimes t^{-x}\ptl_\ell
-(x_j+\d_{j,\ell})t^x\ptl_k\otimes t^{-x+\es_j}\ptl_\ell.
\eqno(3.14)$$
Thus by replacing $D$ by $D-D'$, where $D'$ is an inner derivation which is a
combination of some $(t^x\ptl_k\otimes t^{-x}\ptl_\ell)_d$ with
$x$ satisfying (3.13) (hence after this replacement, (3.11) still holds),
we can suppose
that if $c^{(1)}_{x,k,\ell}\ne0$ then
$$
\left\{
\begin{array}{lll}
\mbox{(i) }
x_j\notin\Z\mbox{ \ and \ }
c^{(1)}_{x+m\es_j,k,\ell}=0
\mbox{ \ for all }m\in\Z\bs\{0\},\mbox{ \ or}\\[0pt]
\mbox{(ii) }x_j=0,-1\mbox{ \ and \ }k\ne j\ne\ell,\mbox{ \ or}\\[0pt]
\mbox{(iii) }x_j=0,-2\mbox{ \ and \ }k\ne j=\ell,\mbox{ \ or}\\[0pt]
\mbox{(iv) }x_j=1,-1\mbox{ \ and \ }k=j\ne\ell,\mbox{ \ or}\\[0pt]
\mbox{(v) }x_j=1,-2\mbox{ \ and \ }k=j=\ell.
\end{array}
\right.
\eqno(3.15)$$
Applying $D$ to
$[t^{-\es_j}\ptl_j,t^{\es_j}\ptl_j]=2\ptl_j$,
and
$[t^{-2\es_j}\ptl_j,t^{\es_j}\ptl_j]=3t^{-\es_j}\ptl_j$,
we obtain
$$
\begin{array}{ll}
(x_j+1+\d_{j,k})c^{(1)}_{x,k,\ell}
\!\!\!\!&-(x_j-1-\d_{j,\ell})c^{(1)}_{x-\es_j,k,\ell}
\\[2pt]&
-(x_j-1-\d_{j,k})c^{(-1)}_{x,k,\ell}
+(x_j+1+\d_{j,\ell})c^{(-1)}_{x+\es_j,k,\ell}=0,
\\[4pt]
(x_j+1+2\d_{j,k})c^{(1)}_{x,k,\ell}
\!\!\!&-(x_j-2-2\d_{j,\ell})c^{(1)}_{x-2\es_j,k,\ell}
\\[2pt]&
-(x_j-2-\d_{j,k})c^{(-2)}_{x,k,\ell}
+(x_j+1+\d_{j,\ell})c^{(-2)}_{x+\es_j,k,\ell}=
3c^{(-1)}_{x,k,\ell\,}.
\end{array}
\eqno\begin{array}{r}\ \\[2pt](3.16)\\[4pt]\ \\[2pt](3.17)
\end{array}
$$
We shall prove
$$
c^{(1)}_{x,k,\ell}=
c^{(-1)}_{x,k,\ell}=c^{(-2)}_{x,k,\ell}=0\mbox{ \ \ for all \ }x\in A,\,k,\ell\in\Z.
\eqno(3.18)$$
If $c^{(1)}_{x,k,\ell}=0$ for all triples
$(x,k,\ell)$, then (3.18) will be deducted from (3.16)--(3.17)
as can be seen in the following proof. Thus assume that
$$
a:=c^{(1)}_{x,k,\ell}\ne0
\mbox{ \ \ for some \ }(x,k,\ell)\in A\times\Z^2.
\eqno(3.19)$$
Now we consider the following cases.
\vskip4pt

\ul{\it Case 1: $x_j\notin\Z$.}
Using (3.15) (i) and
the fact that $c^{(-1)}_{x+m_0\es_j,k,\ell}\ne0$ for all but a finite
number of $m$, by replacing $x$ by $x+m\es_j$ in (3.16) and
letting $m=2,3,...$ or $m=-1,-2,,...,$
we obtain
$$
c^{(-1)}_{x+m\es_j,k,\ell}=0\mbox{ \ \ for \ }m\ne1.
\eqno(3.20)$$
Then by replacing $x$ by $x+m\es_j$ for $m=0,1$ in (3.16), we obtain
$$
\left\{
\begin{array}{rrr}
(x_j+1+\d_{j,k})c^{(1)}_{x,k,\ell}
\!\!\!&
+(x_j+1+\d_{j,\ell})c^{(-1)}_{x+\es_j,k,\ell}
\!\!\!&
=0,
\\[2pt]
-(x_j-\d_{j,\ell})c^{(1)}_{x,k,\ell}
\!\!\!&
-(x_j-\d_{j,k})c^{(-1)}_{x+\es_j,k,\ell}
\!\!\!&=0.
\end{array}\right.
\eqno(3.21)$$
This together with (3.19) shows that the determinant of coefficients of
$c^{(1)}_{x,k,\ell}$ and
$c^{(-1)}_{x+\es_j,k,\ell}$ in (3.21) is zero, which gives $\d_{j,k}=\d_{j,\ell}$, and
so $c^{(-1)}_{x+\es_j,k,\ell}=-a$. Thus there are two subcases:
\vskip4pt

\ul{\it Subcase (1.i): $j=k=\ell$.}
Similar to the proof of (3.20), by replacing $x$ by $x+m\es_j$ in (3.17) and letting $m=3,4,...$ or $m=-1,-2,...$,
we obtain $c^{(-2)}_{x+m\es_j,k,\ell}=0$ if $m\ne1,2$, and
then by letting $m=0,1,2$ and using (3.19), we obtain
$$\left\{
\begin{array}{ll}
(x_j+3) a
+
(x_j+2)c^{(-2)}_{x+\es_j,k,\ell}=0,
\\[2pt]
-(x_j-2)a-(x_j-1)c^{(-2)}_{x+2\es_j,k,\ell}=0,
\\[2pt]
-(x_j-2)c^{(-2)}_{x+\es_j,k,\ell}
+(x_j+3)c^{(-2)}_{x+2\es_j,k,\ell}=-3a.
\end{array}\right.
\eqno(3.22)$$
Regarding (3.22) as a system of linear equations on variables
$a,\,c^{(-2)}_{x+\es_j,k,\ell}$ and
$c^{(-2)}_{x+2\es_j,k,\ell{\ssc\,}}$,
it is immediate to compute the determinant of coefficients is nonzero, this
is a contradiction with $a\ne0$. Thus this subcase does not occur.
\vskip4pt

\ul{\it Subcase (1.ii): $k\ne j\ne\ell$.}
Since $M_u$ is a finite set for any $u\in\WW$, we can
fix $s>>0$ such that
$$
x_k-s\ne0\mbox{ \ \ and \ }
c^{(1)}_{x+m\es_j+s\es_k,k,\ell}=0\mbox{ \ \ for all \ }m\in\Z.
\eqno(3.23)$$
Recall (3.12). When $u=t^{s\es_k}\ptl_k$, we shall
simply denote the coefficient $c^{(u)}_{x',k',\ell'}$
by $b_{x',k',\ell'}$ for $(x',k',\ell')\in M_u$.
Applying $D$ to
$$
[t^{s\es_k}\ptl_k,t^{\es_j}\ptl_j]=0,
\eqno(3.24)$$
and comparing
the coefficients of $t^{x+(m+1)\es_j+s\es_k}\ptl_k\otimes
t^{-x-m\es_j}\ptl_{\ell\ssc\,}$, using (3.23), we obtain
$$
(x_k-s)c^{(1)}_{x+m\es_j,k,\ell}
=(x_j+m)(b_{x+m\es_j,k,\ell}-b_{x+m\es_j+\es_j,k,\ell}).
\eqno(3.25)$$
As
the proof of (3.20), by letting $m=1,2,...$ or $m=-1,-2,...$,
we obtain $b_{x+m\es_j,k,\ell}=0$ for $m\in\Z$. Then (3.25)
shows that $c^{(1)}_{x,k,\ell}=0$, contradicting (3.19).
Thus this subcase does not occur.
\vskip4pt

\ul{\it Case 2: $x_j\in\Z$.} Replacing $x$ by $x-x_j\es_j$ we can
suppose $x_j=0$ such that there exist at
most two $m$'s with $c^{(1)}_{x+m\es_j,k,\ell}\ne0$ (see (3.15)).
We consider the following subcases.
\vskip4pt

\ul{\it Subcase (2.i): $k\ne j\ne\ell$.}
Then (3.15) means
$$
c^{(1)}_{x+m\es_j,k,\ell}=0\mbox{ \ \  if \ }m\ne0,-1.
\eqno(3.26)$$
As in Subcase (1.ii), we still have (3.25), which in turn shows that
$$
b_{x+m\es_j,k,\ell}=0\mbox{ if }m\ne0
\mbox{ \ (in particular $c^{(1)}_{x,k,\ell}=0$), and \ }
c^{(1)}_{x-\es_j,k,\ell}=(x_k-s)^{-1}b_{x,k,\ell}.
\eqno(3.27)
$$
Replacing $t^{\es_j}\ptl_j$ by $t^{-2\es_j}\ptl_j$ in (3.24), then
(3.25) becomes (recall that we suppose $x_j=0$)
$$
(x_k-s)c^{(-2)}_{x+m\es_j,k,\ell}
=m(b_{x+m\es_j,k,\ell}-b_{x+m\es_j-2\es_j,k,\ell}).
\eqno(3.28)$$
Now (3.28) together with (3.27) shows that
$$
c^{(-2)}_{x+m\es_j,k,\ell}=0\mbox{ \ if \ }m\ne2,\mbox{ \ and \ }
c^{(-2)}_{x+2\es_j,k,\ell}=-2(x_k-s)^{-1}b_{x,k,\ell}.
\eqno(3.29)$$
Similarly, we have $c^{(-1)}_{x+m\es_j,k,\ell}=0$ if $m\ne1$ and
$c^{(-1)}_{x+\es_j,k,\ell}=-(x_k-s)^{-1}b_{x,k,\ell}$. This together
with (3.29) and (3.17) proves that $b_{x,k,\ell}=0$, and thus
(3.18) holds.
\vskip4pt

\ul{\it Subcase (2.ii): $\d_{j,k}\ne\d_{j,\ell}$.}
By symmetry, we can suppose $k\ne j=\ell$. Then (3.15) means that
$c^{(1)}_{x+m\es_j,k,j}=0$ if $m\ne0,-2$. In this case,
(3.25) becomes
$$
(x_k-s)c^{(1)}_{x+m\es_j,k,j}
=mb_{x+m\es_j,k,j}-(m+1)b_{x+m\es_j+\es_j,k,j}.
\eqno(3.30)$$
Setting $m=1,2,...$ gives $b_{x+m\es_j,k,j}=0$ if $m\ge1$.
Setting
$m=0$ gives $c^{(1)}_{x,k,j}=0$. Setting $m=-1$ gives
$b_{x-\es_j,k,j}=0$. Setting $m=-3,-4,...$ gives
$b_{x+m\es_j,k,j}=0$ if $m\le-2$. Finally setting $m=-2$ gives
$c^{(1)}_{x-2\es_j,k,j}=0$.
This proves $c^{(1)}_{x+m\es_j,k,j}=0$ for all $m\in\Z$, thus (3.18) can be deducted.
\vskip4pt

\ul{\it Subcase (2.iii): $j=k=\ell$.}
For simplicity, we denote $c^{(p)}_m=c^{(p)}_{x+m\es_j,j,j}$.
Then (3.15) means that
$$
c^{(1)}_m=0\mbox{ \ \ if \ }m\ne1,-2.
\eqno(3.31)$$
From this and (3.16), we obtain
$$
c^{(-1)}_m=0\mbox{ \ if \ }m\ne2,-1
\mbox{ \ and \ }
c^{(-1)}_m=-c^{(1)}_{m-1}
\mbox{ \ if \ }m=2,-1.
\eqno(3.32)$$
Then (3.17) becomes
$$
(m+3)c^{(1)}_m-(m-4)c^{(1)}_{m-2}-(m-3)c^{(-2)}_m+(m+2)c^{(-2)}_{m+1}
=3c^{(-1)}_m=-3c^{(1)}_{m-1}.
\eqno(3.33)$$
Setting $m=4,5,...$ gives $c^{(-2)}_m=0$ if $m\ge 4$. Setting
$m=-3,-4,...$ gives $c^{(-2)}_m=0$ if $m\le -2$. Setting $m=3,-2$ gives
$c^{(1)}_1=c^{(1)}_{-2}=0$. Setting $m=2,1,0,-1$ gives
$$
c^{(-2)}_2+4c^{(-2)}_3=0,\ \ 2c^{(-2)}_1+3c^{(-2)}_2=0,\ \
3c^{(-2)}_0+2c^{(-2)}_1=0,\ \ 4c^{(-2)}_{-1}+c^{(-2)}_0=0.
\eqno(3.34)$$
If we take
$$
u=t^{x+\es_j}\ptl_j\otimes t^{-x-\es_j}-2
t^x\ptl_j\otimes t^{-x}+t^{x-\es_j}\ptl_j\otimes t^{-x+\es_j},
\eqno(3.35)$$
then it is straightforward to verify $u_d(t^{\pm\es_j}\ptl_j)=0$.
Thus if we replace $D$ by $D-c^{(-2)}_3u_d$ (cf.~the statement
after (3.14)), we can suppose $c^{(-2)}_3=0$.
This together with the above results shows that we have (3.18).
This completes the proof of Claim 2.
\vskip4pt

{\bf Claim 3}. $D=0$.

First from Claim 2, we can easily
deduct that $D(t^{p{\ssc\,}\es_j}\ptl_j)=0$ for
$j=1,2,...,n$ and $p\in\Z$. For any $u=t^z\ptl_j\in\WW_z$,
assume we have (3.12). Fix $p_0>>0$ such that
$$
c^{(u)}_{x\pm p{\ssc\,}\es_j,k,\ell}=0\mbox{ \ \ for all \ }
p>p_0,\,x\in A,\,
k,\ell=1,2,...,n.
\eqno(3.36)$$
Applying $D$ to
$$
[t^{-p{\ssc\,}\es_j}\ptl_j,[t^{p{\ssc\,}\es_j}\ptl_j,t^z\ptl_j]]
=(z_j-p)(z_j+2p)t^z\ptl_j,
$$
and comparing the coefficients of $t^{x+z}\ptl_k\otimes
t^{-x}\ptl_\ell{\ssc\,}$, using (3.36), we
obtain
$$
((x_j+z_j-p)(x_j+z_j+2p)+(-x_j-p)(-x_j+2p)-(z_j-p)(z_j-2p))
c^{(u)}_{x,k,\ell}=0 \mbox{ for all }p>p_0.
$$
Thus $c^{(u)}_{x,k,\ell}=0$. This proves the claim and the
theorem. \hfill$\Box$ \vskip7pt \ni{\bf Lemma 3.7.} {\it Suppose
$r\in V$ such that $a\cdot r\in {\rm Im}(1-\tau)$ for all
$a\in\WW$. Then $r\in{\rm Im}(1-\tau)$.} \vskip5pt

\ni{\it Proof.} First note that $\WW\cdot {\rm
Im}(1-\tau)\subset{\rm Im}(1-\tau)$. We shall prove that after a
number of steps in each of which $r$ is replaced by $r-u$ for some
$u\in{\rm Im}(1-\tau)$, the zero element is obtained and thus proving that
$r\in{\rm Im}(1-\tau)$.

 Write $r=\sum_{x\in A}r_x$. Obviously,
$$ r\in{\rm Im}(1-\tau)\ \ \Longleftrightarrow\ \ r_x\in{\rm
Im}(1-\tau)\mbox{ \ \ for all \ }x\in A. \eqno(3.37)$$ For any
$x'\ne0$, choose $\ptl\in T$ such that $\ptl(x')\ne0$. Then
$\sum_{x\in A}\ptl(x)r_x=\ptl\cdot r\in{\rm Im}(1-\tau)$. By
(3.37), $\ptl(x)r_x\in{\rm Im}(1-\tau)$, in particular,
$r_{x'}\in{\rm Im}(1-\tau)$. Thus by replacing $r$ by
$r-\sum_{0\ne x\in A}r_x$, we can suppose $r=r_0\in V_0$.

Now write $r=\sum c_{x,k,\ell{\ssc\,}}t^x\ptl_k\otimes t^{-x}\ptl_\ell$.
Choose a total order on $A$ compatible with its group structure.
Since $a_{x,k,\ell}:=t^x\ptl_k\otimes
t^{-x}\ptl_\ell-t^{-x}\ptl_\ell\otimes t^{x}\ptl_k\in{\rm
Im}(1-\tau)$, by replacing $r$ by $r-u$, where $u$ is some
combination of $a_{x,k,\ell}$, we can suppose
$$
c_{x,k,\ell}\ne0\ \ \Longrightarrow\ \ x\ge0\mbox{ \ or \ }x=0,\
k\le \ell. \eqno(3.38)$$ Now assume that $c_{x,k,\ell}\ne0$ for
some $x>0$. Fix $m>>0$ if $\es_k>0$, or $m<<0$ if $\es_k<0$, such
that $x_k-m\ne0$, then we see that the term
$t^{x+m\es_k}\ptl_k\otimes t^{-x}\ptl_\ell$ appears in $
t^{m\es_k}\ptl_k\cdot r$, but (3.38) implies that
the term $t^{-x}\ptl_\ell\otimes
t^{x+m\es_k}\ptl_k$ does not appear in $ t^{m\es_k}\ptl_k\cdot r$,
a contradiction with the fact that $ t^{m\es_k}\ptl_k\cdot
r\in{\rm Im}(1-\tau)$. Similarly, we can obtain a contradiction if
$c_{0,k,\ell}\ne0$ for some $k\le\ell$. Thus by (3.38), $r=0$.
\hfill$\Box$
 \vskip7pt

 Now we can obtain the main result of this paper.
\vskip7pt

\ni{\bf Theorem 3.8.} {\it Let $(\WW,[\cdot,\cdot])$ be the Lie
algebra of generalized Witt type. Then every Lie bialgebra
structure on $\WW$ is a triangular coboundary Lie bialgebra. }
\vskip5pt

\ni{\it Proof.} Let $(\WW,[\cdot,\cdot],\D)$ be a Lie bialgebra
structure on $\WW$. By Definition 2.3(iii) and Theorem 3.6,
$\D=\D_r$ is defined by (\ref{add2.0}) for some
$r\in\WW\otimes\WW$. By Definition 2.2(i), ${\rm
Im}\,\D\subset{\rm Im}(1-\tau)$. Thus by Lemma 3.7, $r\in{\rm
Im}(1-\tau)$, namely, $r=\sum_{i=1}^m (a_i\otimes b_i-b_i\otimes
a_i)$ for some $a_i,b_i\in\WW$. Then Definition 2.2(ii),
(\ref{add-c}) and Proposition 3.5 show that $c(r)=0$. Thus
Definitions 2.5 and 2.6 say that $(\WW,[\cdot,\cdot],\D)$ is a
triangular coboundary Lie bialgebra.
\hfill$\Box$
\vskip12pt

\cl{\bf References}\vs{0pt}

\vskip5pt\small
\parindent=8ex\parskip=2pt\baselineskip=2pt
\re{D} V.G. Drinfeld, Quantum groups, in: {\it Proceeding of the
International Congress of Mathematicians}, Vol.~1, 2, Berkeley,
Calif.~1986, Amer.~Math.~Soc., Providence, RI, 1987, pp.~798-820.

 \re{DZ} D. Dokovic and K. Zhao, Derivations, isomorphisms
and second cohomology of generalized Witt algebras, {\it Trans.
Amer. Math. Soc.} {\bf350} (1998), 643--664.

\re{F} R. Farnsteiner, Derivations and central extensions of
finitely generalized Lie algebras, {\it J. Algebra} {\bf 118}
(1988), 33--45.

\re{M} %\re{M1} 
W. Michaelis, A class of infinite-dimensional Lie
bialgebras containing the Virasoro algebras, {\it Adv. Math.}
{\bf107} (1994), 365--392.

%\re{M2} W. Michaelis, Lie coalgebras, {\it Adv. Math.} {\bf38}
%(1980), 1--54.
%
%\re{M3} W. Michaelis, The dual Poincare-Birkhoff-Witt theorem,
%{\it Adv. Math.} {\bf57} (1985), 93--162.

\re{NT} S.-H. Ng and E.J. Taft, Classification of the Lie
bialgebra structures on the Witt and Virasoro algebras,
{\it J. Pure Appl.~Algebra} {\bf151} (2000), 67--88.

\re{N} W.D. Nichols, The structure of the dual Lie coalgebra
of the Witt algebra, {\it J. Pure Appl. Algebra} {\bf 68} (1990),
395--364.

\re{P}  D. Passman,  New simple infinite dimensional Lie
algebras, {\it J. Algebra} {\bf 206} (1998), 682--692.

\re{SXZ} Y. Su, X. Xu and H. Zhang, Derivation-simple
algebras and the structures of  Lie algebras of Witt type,
{\it J. Algebra} {\bf 233} (2000), 642--662.

\re{SZ} Y. Su and K. Zhao, The second cohomology group of
generalized of Witt type Lie algebras and certain representations,
{\it Comm. Algebra} {\bf30} (7) (2002), 3285--3309.

\re{T} E.J. Taft, Witt and Virasoro algebras as Lie
bialgebras, {\it J. Pure Appl.~Algebra} {\bf87} (1993), 301--312.

\re{X} X. Xu, New generalized simple Lie algebras of Cartan
type over a field with characteristic 0,  {\it J.~Algebra} {\bf 224}
(2000), 23--58.

\end{document}